\documentclass[10pt]{amsart}

\newcommand{\Fq}{{\mathbf F}_q}
\newcommand{\Fqq}{{\mathbf F}_{q^2}}
\newcommand{\BP}{{\mathbf P}}
\newcommand{\BQ}{{\mathbf Q}}
\newcommand{\BF}{{\mathbf F}}
\newcommand{\BZ}{{\mathbf Z}}

\newcommand{\Kplus}{K^{+}}
\DeclareMathOperator{\End}{End}

\newtheorem*{theoremstar}{Theorem}

\begin{document}

\title{Appendix to a paper of Maisner and Nart}
\author{Everett W.\ Howe}
\address{Center for Communications Research, 
         4320 Westerra Court, 
         San Diego, CA 92121-1967, USA.}
\email{however@alumni.caltech.edu}
\urladdr{http://alumni.caltech.edu/\~{}however/}
\date{31 October 2001}


\begin{abstract}
We prove that
there is no genus-$2$ curve over $\Fq$ whose Jacobian has characteristic
polynomial of Frobenius equal to $x^4 + (1-2q)x^2 + q^2$.
\end{abstract}

\maketitle
For every prime power $q$ let $f_q$ denote the polynomial $x^4 + (1-2q)x^2 + q^2$.
In~\cite{mn},
Maisner and Nart observe that for all prime 
powers $q\le 64$, no genus-$2$ curve over $\Fq$ has characteristic polynomial~$f_q$.
(By the {\em characteristic polynomial\/} of a curve we mean the characteristic
polynomial of the Frobenius endomorphism of the Jacobian of
the curve.)  The purpose of this appendix is to prove that Maisner and Nart's
observation holds for all prime powers~$q$.

\begin{theoremstar}
\label{T-main}
There is no curve of genus $2$ over any finite field $\Fq$ whose characteristic
polynomial is equal to~$f_q$.
\end{theoremstar}

\begin{proof}
Suppose, to obtain a contradiction, that $C$ is a genus-$2$ curve over 
a finite field $\Fq$
whose characteristic polynomial is equal to~$f_q$. Note that then
$\#C(\Fq) = q + 1$ and $\#C(\Fqq) = (q-1)(q-3)$.

Let $J$ be the Jacobian of $C$, let $\lambda$ be the canonical principal
polarization of~$J$, let $F$ be the Frobenius endomorphism of~$J$,
and let $V = q/F$ be the Verschiebung endomorphism of~$J$.  Since $f_q$ is
irreducible and its middle coefficient is coprime to $q$, we see that $J$ is
a simple ordinary abelian surface, and it follows that $K = (\End J)\otimes\BQ$
is equal to the field~$\BQ(F)$.  In fact, $K$ is a totally imaginary quadratic 
extension of a totally real quadratic field $\Kplus$, and general theory
(see~\cite[p.~201]{mumford}) shows that the Rosati involution
$x \mapsto x^\dagger$ on $K$ is complex conjugation.

Let $i$ be the endomorphism $F-V$ of $J$.  It is easy to check that $i^2 = -1$,
and it follows that $i^\dagger i = 1$.  Thus $i$ is an automorphism of $J$
that respects the polarization~$\lambda$, so $i$ can be viewed as an automorphism of
the polarized abelian variety~$(J,\lambda)$.
Since $C$ is hyperelliptic, Torelli's theorem (see~\cite[p.~202]{milne}) shows
that the natural map from the automorphism group of $C$ to the automorphism
group of $(J,\lambda)$ is an isomorphism that takes the hyperelliptic involution
to~$-1$.  Thus, the automorphism $i$ of $(J,\lambda)$ gives us
an automorphism $\alpha$ of $C$, defined over~$\Fq$,
whose square is the hyperelliptic involution.

Let $W$ denote the set of Weierstrass points of $C$, viewed as a set with an
action of the absolute Galois group of $\Fq$. 
If $P$ is a geometric point of $C$ whose orbit under the action of
$\alpha$ contains fewer than four points, then $P$ must be fixed
by $\alpha^2 = -1$, so $P$ must lie in $W$.
Thus, for every finite extension field $k$ of $\Fq$ we have 
$\#C(k) \equiv \#W(k) \bmod 4.$

Suppose that $q$ is odd.  Then $W$ consists of six points, and we will show
that exactly two of these points are fixed by $\alpha$.

Consider the map $C\to\BP^1$ obtained from the hyperelliptic involution,
and let $W'$ denote the set of six points of $\BP^1$ lying under the 
Weierstrass points of $C$.  The automorphism $\alpha$ induces an involution
$\beta$ of $\BP^1$ that takes the set $W'$ to itself.
Geometrically, this involution is conjugate to the involution $x\mapsto -x$,
so if none of the points in $W'$ were fixed by $\beta$ the curve $C$
would be isomorphic (over the algebraic closure of $\Fq$)
to a curve of the form $y^2 = f(x^2)$, where $f$ is a cubic polynomial.
But then $\alpha$ would have to be of the form $(x,y)\mapsto(-x,\pm y)$,
and such an automorphism has order two.  Thus, $\beta$ must fix at least
one of the six points of~$W'$.  But the points not fixed by $\beta$ come
in plus/minus pairs, so there must be at least two points of $W'$ fixed
by~$\beta$.  Since $x\mapsto-x$ has exactly two fixed points in~$\BP^1$, 
there must be exactly two points of $W'$ fixed by $\beta$.  
It follows that exactly two points of $W$ are fixed by $\alpha$, as claimed.

Since $\alpha$ is defined over $\Fq$, the two points of $W$ fixed by $\alpha$
must be defined over~$\Fqq$.  Thus, $\#W(\Fqq) \ge 2$.  But we also have
$$\#W(\Fqq) \equiv \#C(\Fqq) = (q-1)(q-3) \equiv 0 \bmod 4,$$
so we must have $\#W(\Fqq) = 4$.  But this is impossible, as one can see
by asking where the other two points of $W$ are defined.
This proves the theorem for odd $q$.

Suppose that $q$ is a power of $2$.  Then $q$ must be a multiple of~$4$, because 
if $q$ were $2$ the curve $C$ would have $-1$ points over $\BF_4$.  
We see that $\#W(\Fq) \equiv 1 \bmod 4$ and $\#W(\Fqq) \equiv 3 \bmod 4$.
But a genus-$2$ curve in characteristic $2$ has at most three Weierstrass points,
so $C$ must have exactly three Weierstrass points, and exactly one of them is
defined over~$\Fq$.  

Once again we let $W'$ denote the points of $\BP^1$ lying under the
Weierstrass points of $C$ and we let $\beta$ be the involution of $\BP^1$
obtained from~$\alpha$.
Clearly $\beta$ must fix the unique point of $W'(\Fq)$.
But $\beta$ cannot fix the other two points of $W'$, because in that case
$\beta$ would be the identity on $\BP^1$, and $\alpha$ could not have
order four.  Thus, $\beta$ must swap the other two points of~$W'$.
It follows that over the algebraic closure of $\Fq$ we can write $C$ as
$y^2 + y = ax + b/x + b/(x+1)$, where we have chosen the coordinates so
that $\beta$ is given by $x\mapsto x+1$.  But then $\alpha$ must send 
$(x,y)$ to $(x+1,y+c)$ where $c^2 + c = a$, and this automorphism has 
order two.  Once again we obtain a contradiction, and the theorem is proved
for even $q$ as well.
\end{proof}

Maisner and Nart also note that for every odd prime power $q<64$,
no genus-$2$ curve over $\Fq$ has characteristic polynomial
$g_q = x^4 + (2-2q)x^2 + q^2$.
The obvious conjecture
is that the same is true for all odd prime powers~$q$.
Unfortunately, the argument we used above cannot be easily modified 
to prove this conjecture; the critical fact we used was that the ring
$\BZ[F,V]$ contains a root of unity other than~$\pm1$, and this is no longer
true when we replace $f_q$ with $g_q$ in our argument.

\end{document}